\RequirePackage{fix-cm}
\documentclass[smallcondensed]{svjour3}     
\smartqed  
\usepackage[parfill]{parskip}    
\usepackage{amsmath, amssymb, latexsym, enumerate, mathrsfs}
\usepackage{graphicx}
\usepackage{graphicx}
\usepackage{array}
\usepackage{times}
\usepackage{moreverb}
\usepackage{marvosym}
\DeclareMathOperator*{\argmin}{arg\,min}
\def\qed{\hfill $\Box$}

\usepackage{latexsym}
\journalname{}
\begin{document}

\title{On noncooperative $n$-player principal eigenvalue games
}

\titlerunning{On noncooperative $n$-player principal eigenvalue games}        

\author{Getachew K. Befekadu \and Panos~J.~Antsaklis
}

\institute{G. K. Befekadu~ ({\large\Letter}\negthinspace) \at
	          Department of Electrical Engineering, University of Notre Dame, Notre Dame, IN 46556, USA. \\
	          Tel.: +1 574 631 6618\\
	          Fax: +1 574 631 4393 \\
	          \email{gbefekadu1@nd.edu}           
	          \and
	          P. J. Antsaklis \at
	          Department of Electrical Engineering, University of Notre Dame, Notre Dame, IN 46556, USA. \\
	          \email{antsaklis.1@nd.edu}           
}

\date{Received: May 21, 2014 / Accepted: date}

\maketitle

\begin{abstract}
We consider a noncooperative $n$-player principal eigenvalue game which is associated with an infinitesimal generator of a stochastically perturbed multi-channel dynamical system -- where, in the course of such a game, each player attempts to minimize the asymptotic rate with which the controlled state trajectory of the system exits from a given bounded open domain. In particular, we show the existence of a Nash-equilibrium point (i.e., an $n$-tuple of equilibrium linear feedback operators) that is distinctly related to a unique maximum closed invariant set of the corresponding deterministic multi-channel dynamical system, when the latter is composed with this $n$-tuple of equilibrium linear feedback operators.
\end{abstract}

\section{Introduction}	 \label{S1}
In this paper, we consider a noncooperative $n$-player principal eigenvalue game which is associated with an infinitesimal generator pertaining to the following stochastically perturbed multi-channel dynamical system\footnote{e.g., see \cite{BefGA13b} for additional discussion on multi-channel dynamical systems without random perturbation terms.}
\begin{align}
d x^{\epsilon}(t) = A x^{\epsilon}(t) dt + \sum\nolimits_{i=1}^n B_i u_i(t) dt + \sqrt{\epsilon} \sigma(x^{\epsilon}(t))dW(t), \,\, x^{\epsilon}(0) = x_0, \label{Eq1}
\end{align}
where
\begin{itemize}
\item[-] $A \in \mathbb{R}^{d \times d}$, $B_i \in \mathbb{R}^{d \times r_i}$, $\epsilon$ is a small positive number (which represents the level of random perturbation in the system),
\item[-] $\sigma \colon \mathbb{R}^{d} \rightarrow \mathbb{R}^{d \times d}$ is Lipschitz with the least eigenvalue of $\sigma(\cdot)\sigma^T(\cdot)$ uniformly bounded away from zero, i.e., 
\begin{align*}
 \sigma(x)\sigma^T(x)  \ge \kappa I_{d \times d} , \quad \forall x \in \mathbb{R}^{d},
\end{align*}
for some $\kappa > 0$,
\item[-] $W(\cdot)$ is a $d$-dimensional standard Wiener process,
\item[-] $x^{\epsilon}(\cdot) \in \mathcal{X} \subseteq \mathbb{R}^{d}$ is the state trajectory of the system,
\item[-] $u_i(\cdot)$ is a $\,\mathcal{U}_i$-valued measurable control process to the $i$th-channel (i.e., an admissible control from the measurable set $\mathcal{U}_i\subset \mathbb{R}^{r_i}$) such that for all $t > s$, $W(t)-W(s)$ is independent of $u_i(\nu)$ for $\nu \le s$ and 
\begin{align*}
 \int_{0}^{t_1} \vert u_i(t)\vert^2 dt < \infty, \quad \forall t_1 \ge 0,
\end{align*}
for $i=1, 2, \ldots, n$.
\end{itemize}
 
Let $D \subset \mathbb{R}^{d}$ be a bounded open domain with smooth boundary (i.e., $\partial D$ is a manifold of class $C^2$). Moreover, denote by $C_{0T}([0,T], \mathbb{R}^d)$ the space of all continuous functions $\varphi(t)$, $t \in [0,\, T]$, with range in $\mathbb{R}^d$; and, in this space, we define the following metric  
\begin{align}
\rho_{0T}(\varphi, \psi) =\sup_{t \in [0,\, T]} \Bigl\vert \varphi(t) - \psi(t) \Bigr\vert, \label{Eq2}
\end{align}
when $\varphi(t)$, $\psi(t)$ belong to $C_{0T}([0,T], \mathbb{R}^d)$.  If $\Phi$ is a subset of the space $C_{0T}([0,T], \mathbb{R}^d)$, then we define
\begin{align}
d_{0T}(\psi, \Phi) =\sup_{\varphi \in \Phi} \rho_{0T}(\psi(t),\varphi(t)). \label{Eq3}
\end{align}

In what follows, we consider a particular class of admissible controls $u_i(\cdot) \in \mathcal{U}_i$ of the form $u_i(t)=\bigl(\mathcal{K}_i x^{\epsilon}\bigr)(t)$, $\forall t \ge 0$, where $\mathcal{K}_i$, for $i =1, 2, \ldots, n$, is a real, continuous $r_i \times d$ matrix function such that
\begin{align}
  \mathscr{K} \subseteq \Biggl\{\underbrace{\bigl(\mathcal{K}_1, \mathcal{K}_2, \ldots, \mathcal{K}_n\bigr)}_{\substack{\triangleq \mathcal{K}}} \in \prod\nolimits_{i=1}^n \mathscr{K}_i[\mathcal{X}, \mathcal{U}_i] \biggm\lvert \phi\bigl(t; 0, x_0, (\mathcal{K}x^{0})(t)\bigr) \in \Omega, \notag \\ \forall t \ge 0, \, \, \forall x_0 \in \Omega \Biggr\}, \label{Eq4} 
\end{align}
where $\mathscr{K}_i[\mathcal{X}, \mathcal{U}_i]$ is a closed subspace of bounded linear feedback operators from $\mathcal{X}$ to $\,\mathcal{U}_i$, $\Omega$ is a bounded open set which is contained in $D \cup \partial D$; and $\phi\bigl(t; 0, x_0, (\mathcal{K} x^{0})(t)\bigr)$ is the unique solution for
\begin{align}
  \dot{x}^{0}(t)=Ax^{0}(t) + \sum\nolimits_{i=1}^n B_i \bigl(\mathcal{K}_i x^{0}\bigr)(t),  \,\, x^{0}(0) = x_0 \in \Omega,  \label{Eq5} 
\end{align}
that corresponds to the deterministic multi-channel dynamical system, when $\epsilon = 0$.

Further, the infinitesimal generator pertaining to the diffusion process $x^{\epsilon}(t)$ of Equation~\eqref{Eq1}, with $u_i(t)=\bigl(\mathcal{K}_i x^{\epsilon}\bigr)(t)$, for $t \ge 0$ and $i =1, 2, \ldots, n$, is given by
\begin{align}
 \mathcal{L}_{\epsilon}^{\mathcal{K}}(\cdot)(x) = \Bigl \langle \bigtriangledown(\cdot), \Bigl(A x + \bigl(B, \mathcal{K}\bigr)x \Bigr) \Bigr\rangle + \frac{\epsilon}{2} \operatorname{tr}\Bigl \{\sigma(x)\sigma^T(x)\bigtriangledown^2(\cdot) \Bigr\}, \label{Eq6}
\end{align}
where $\bigl(B, \mathcal{K}\bigr)x(\cdot) = \sum\nolimits_{i=1}^n B_i \bigl(\mathcal{K}_i x\bigr)(\cdot)$ for all $t \ge 0$.

For any fixed $\mathcal{K} \in \mathscr{K}$ and $\epsilon > 0$, let $\tau_D^{\epsilon}$ be the first exit-time from the domain $D$ for the diffusion process $x^{\epsilon}(t)$, i.e., 
\begin{align}
\tau_D^{\epsilon} = \inf \bigl\{ t > 0 \, \bigl\vert \, x^{\epsilon}(t) \notin D \bigr\}, \label{Eq7}
\end{align}
which also depends on the class of linear feedback operators $\mathscr{K}$ and, in particular, on the behavior of the solutions to the deterministic dynamical system in Equation~\eqref{Eq5}. Moreover, let us denote by $\lambda_{\epsilon}^{\mathcal{K}}$ the principal eigenvalue of the infinitesimal generator $-\mathcal{L}_{\epsilon}^{\mathcal{K}}$  corresponding to a zero boundary condition on $\partial D$ which is given by
\begin{align}
\lambda_{\epsilon}^{\mathcal{K}} = - \limsup_{T \rightarrow \infty} \frac{1} {T} \log \mathbb{P}_{\epsilon}^{\mathcal{K}} \bigl\{\tau_D^{\epsilon} > T \bigr\}, \label{Eq8}
\end{align}
where the probability $\mathbb{P}_{\epsilon}^{\mathcal{K}}\bigl\{\cdot\bigr\}$ is conditioned on the initial point $x_0 \in D$ as well as on the class of linear feedback operators $\mathscr{K}$.

Next, let us introduce the following definition (i.e., the maximum closed invariant set for the deterministic dynamical system of Equation~\eqref{Eq5} under the action of the class of linear feedback operators $\mathscr{K}$) which is useful in the sequel.
\begin{definition} \label{DFN1}
A set $\Lambda_D^{\mathcal{K}} \subset D \cup \partial D$ is called a maximum closed invariant set (under the action of an $n$-tuple of linear feedback operators with respect to the deterministic dynamical system), if any set $\Omega \subset D \cup \partial D$, for some $\mathcal{K} \in \mathscr{K}$, satisfying the property
\begin{align}
 \phi\bigl(t; 0, x_0, (\mathcal{K} x^{0})(t)\bigr) \in \Omega, \,\, \forall t \ge 0, \,\, \forall x_0 \in \Omega \label{Eq9}
\end{align}
 is a subset of $\Lambda_D^{\mathcal{K}}$.\footnote{Such an invariant set is closed (and it may also be an empty set). Note that the solution $\phi\bigl(t; 0, x_0, (\mathcal{K}x^{0})(t)\bigr)$ corresponds to the deterministic dynamical system, i.e.,
 \begin{align*}
  \dot{x}_t^{0}=Ax^{0}(t) + \sum\nolimits_{i=1}^n B_i \bigl(\mathcal{K}_i x^{0}\bigr)(t),  \,\,  \forall x_0 \in \Omega,  
\end{align*}
when such a system is composed with an $n$-tuple of linear feedback operators $\mathcal{K} \in \mathscr{K}$.}
\end{definition}

In Section~\ref{S2}, we introduce a noncooperative $n$-player principal eigenvalue game -- where, in the course of such a game, each player attempts to minimize a cost criterion related to the asymptotic rate with which the controlled state trajectory of the dynamical system of Equation~\eqref{Eq1} exits from the given bounded open domain $D$. To be specific, we use a cost criterion that is directly related to minimizing the principal eigenvalue of the infinitesimal generator
\begin{align}
\mathscr{K}_i[\mathcal{X}, \mathcal{U}_i] \ni \mathcal{K}_i \mapsto \lambda_{\epsilon, i}^{(\mathcal{K}_i, \mathcal{K}_{\neg i}^{\ast})} \in \mathbb{R}_{+} \cup \{\infty\}, \label{Eq10}
\end{align}
with respect to the admissible controls $u_i(\cdot) \in \mathcal{U}_i$ of the form $u_i(t)=\bigl(\mathcal{K}_i x^{\epsilon}\bigr)(t)$, for $t \ge 0$, where $\mathcal{K}_i  \in \mathscr{K}_i[\mathcal{X}, \mathcal{U}_i]$ for each $i =1, 2, \ldots, n$; while the others $\mathcal{K}_{\neg i}^{\ast}$ remain fixed.\footnote{$(\mathcal{K}_i, \mathcal{K}_{\neg i}^{\ast})\triangleq \bigl(\mathcal{K}_1^{\ast}, \dots, \mathcal{K}_{i-1}^{\ast}, \mathcal{K}_i, \mathcal{K}_{i+1}^{\ast}, \ldots, \mathcal{K}_n^{\ast}\bigr)$.} Note that if the domain $D$ contains an equilibrium point for the deterministic dynamical system of Equation~\eqref{Eq5} (i.e., under the action of the class of linear feedback operators $\mathscr{K}$). Then, the asymptotic behavior of the principal eigenvalue tends to zero exponentially as $\epsilon \rightarrow 0$ (e.g., see \cite{Ven72}, \cite{VenFre70} or \cite{Day83}). On the other hand, if the maximum closed invariant set $\Lambda_D^{(\mathcal{K}_i, \mathcal{K}_{\neg i}^{\ast})} \subset D \cup \partial D$ (under the action of $(\mathcal{K}_i, \mathcal{K}_{\neg i}^{\ast}) \in \mathscr{K}$ with respect to the deterministic dynamical system) is nonempty, then the following asymptotic condition holds 
\begin{align}
-\lim_{\epsilon \rightarrow 0} \, \limsup_{T \rightarrow \infty} \frac{\epsilon} {T} \log \mathbb{P}_{\epsilon, i}^{(\mathcal{K}_i, \mathcal{K}_{\neg i}^{\ast})} \bigl\{\tau_{D, i}^{\epsilon} > T \bigr\} < \infty, \,\,x_0 \in D, \label{Eq11}
\end{align}
where $\tau_{D, i}^{\epsilon}$ the exit time with respect to the $i$th-channel and $(\mathcal{K}_i, \mathcal{K}_{\neg i}^{\ast})\in \mathscr{K}$ for each $i=1,2, \ldots n$. Later, such a connection between the existence of the maximum closed invariant sets $\Lambda_D^{(\mathcal{K}_i, \mathcal{K}_{\neg i}^{\ast})} \subset D \cup \partial D$, for each $i = 1,2, \ldots, n$, for the dynamical system of Equation~\eqref{Eq5} and the asymptotic behavior of the principal eigenvalue of the infinitesimal generator $-\mathcal{L}_{\epsilon}^{\mathcal{K}}$ (which corresponds to a zero boundary condition on $\partial D$) allows us to provide some results on the existence of a Nash-equilibrium point $(\mathcal{K}_1^{\ast}, \mathcal{K}_2^{\ast}, \ldots, \mathcal{K}_n^{\ast}) \in \mathscr{K}$ in a game-theoretic setting (cf. Definition~\ref{DFN2} for the definition of Nash-equilibrium points).

Here, it is worth remaking that the principal eigenvalue $\lambda_{\epsilon, i}^{(\mathcal{K}_i, \mathcal{K}_{\neg i}^{\ast})}$ (for a particular $(\mathcal{K}_i, \mathcal{K}_{\neg i}^{\ast}) \in \mathscr{K}$) is the boundary value between those $R < r_{i}(\mathcal{K}_i, \mathcal{K}_{\neg i}^{\ast})$ for which $\mathbb{E}_{\epsilon, i}^{(\mathcal{K}_i, \mathcal{K}_{\neg i}^{\ast})}\bigl\{\exp(\epsilon^{-1} R \tau_{D, i}^{\epsilon})\bigr\} \\< \infty$ and those $R > r_{i}(\mathcal{K}_i, \mathcal{K}_{\neg i}^{\ast})$ for which $\mathbb{E}_{\epsilon, i}^{(\mathcal{K}_i, \mathcal{K}_{\neg i}^{\ast})}\bigl\{\exp(\epsilon^{-1} R\tau_{D, i}^{\epsilon})\bigr\} = \infty$, where $r_{i}(\mathcal{K}_i, \mathcal{K}_{\neg i}^{\ast})$ is given by the following\footnote{Note that the asymptotic behavior of $\bigl(\epsilon/T\bigr) \log \mathbb{P}_{\epsilon, i}^{(\mathcal{K}_i, \mathcal{K}_{\neg i}^{\ast})} \bigl\{\tau_{D, i}^{\epsilon} > T \bigr\}$, for each $i=1,2, \ldots n$, as $\epsilon \rightarrow 0$ and $T  \rightarrow \infty$, determines whether the deterministic dynamical system in Equation~\eqref{Eq5} has a maximum closed invariant set in $D \cup \partial D$ or not.}
\begin{align}
r_{i}(\mathcal{K}_i, \mathcal{K}_{\neg i}^{\ast}) = \limsup_{T \rightarrow \infty} \,\inf_{\substack{\varphi(t) \in C_{0T}([0,T], \mathbb{R}^d)\\ \varphi(0)=x_0}} \frac{1}{T} \biggl \{ S_{0T}^{(\mathcal{K}_i, \mathcal{K}_{\neg i}^{\ast})}(\varphi(t)) \, \Bigl \vert \, \varphi(t) \in D \cup \partial D,\, t \in [0, T] \biggr\}, \label{Eq12}
\end{align}
with $(\mathcal{K}_i, \mathcal{K}_{\neg i}^{\ast})\in \mathscr{K}$ for each $i=1,2, \ldots n$.

Note that, in general, such an asymptotic analysis involves minimizing the following {\it action functional}
\begin{align}
 S_{0T}^{(\mathcal{K}_i, \mathcal{K}_{\neg i}^{\ast})}(\varphi(t)) =  \frac{1}{2} \int_{0}^{T}\biggl\Vert\frac{d\varphi(t)}{dt} - \Bigl(A \varphi(t) + \bigl(B, (\mathcal{K}_i, \mathcal{K}_{\neg i}^{\ast})\bigr)\varphi(t) \Bigr) \biggr\Vert^2dt, \label{Eq13}
\end{align}
where
\begin{align}
\biggl\Vert\frac{d\varphi(t)}{dt} -  \Bigl(A \varphi(t) + \bigl(B, (\mathcal{K}_i, \mathcal{K}_{\neg i}^{\ast})\bigr)\varphi(t) \Bigr) \biggr\Vert^2 = \biggl[\frac{d\varphi(t)}{dt} - \Bigl(A \varphi(t) + \bigl(B, (\mathcal{K}_i, \mathcal{K}_{\neg i}^{\ast})\bigr)\varphi(t) \Bigr) \biggr]^T \notag \\
 \times \Bigl(\sigma(\varphi(t))\sigma^T(\varphi(t))\Bigr)^{-1} \biggl[\frac{d\varphi(t)}{dt} -  \Bigl(A \varphi(t) + \bigl(B, (\mathcal{K}_i, \mathcal{K}_{\neg i}^{\ast})\bigr)\varphi(t) \Bigr) \biggr], \label{Eq14}
\end{align}
with $\bigl(B, (\mathcal{K}_i, \mathcal{K}_{\neg i}^{\ast})\bigr)\varphi(t) = B_i \bigl(\mathcal{K}_i \varphi\bigr)(t) + \sum\nolimits_{j \neq i} B_j \bigl(\mathcal{K}_j^{\ast} \varphi\bigr)(t)$, for each $i=1,2, \ldots n$, and  where $\varphi(t) \in C_{0T}([0,T], \mathbb{R}^d)$ is absolutely continuous.

In the remainder of this section, we state the following lemmas that will be useful for proving our main results (see \cite[Theorem~1.1, Theorem~1.2 and Lemma~9.1]{VenFre70} or \cite{Ven73}; and see \cite[pp.\,332--340]{Fre76} for additional discussions).
 
\begin{lemma} \label{L1}
For any $\alpha >0$, $\delta > 0$ and $\gamma > 0$, there exists an $\epsilon_0 > 0$ such that
\begin{enumerate} [(i)]
\item
\begin{align}
\mathbb{P}_{\epsilon, i}^{(\mathcal{K}_i, \mathcal{K}_{\neg i}^{\ast})} \Bigl\{\rho_{0T}\bigl(x^{\epsilon}(t), \varphi(t)\bigr) < \delta \Bigr\} \ge \exp \Bigl\{ - \epsilon^{-1}\bigl(S_{0T}^{(\mathcal{K}_i, \mathcal{K}_{\neg i}^{\ast})}(\varphi(t)) + \gamma \bigr)\Bigl\}, \, \forall  \epsilon \in (0, \epsilon_0), \label{Eq15}
\end{align}
where $\varphi(t)$ is any function in $C_{0T}([0,T], \mathbb{R}^d)$ for which $S_{0T}^{(\mathcal{K}_i, \mathcal{K}_{\neg i}^{\ast})}(\varphi(t)) < \alpha$ and $\varphi(0)=x_0$; and
\item
\begin{align}
\mathbb{P}_{\epsilon, i}^{(\mathcal{K}_i, \mathcal{K}_{\neg i}^{\ast})} \Bigl\{\rho_{0T}\bigl(x^{\epsilon}(t), \Phi_{x_0, \alpha}^{(\mathcal{K}_i, \mathcal{K}_{\neg i}^{\ast})} \bigr) \ge \delta \Bigr\} \le \exp \Bigl\{ - \epsilon^{-1}\bigl(\alpha - \gamma \bigr)\Bigl\}, \,\, \forall  \epsilon \in (0, \epsilon_0), \label{Eq16}
\end{align}
where
\begin{align}
\Phi_{x_0, \alpha}^{(\mathcal{K}_i, \mathcal{K}_{\neg i}^{\ast})} = \Bigl\{\varphi(t) \in C_{0T}([0,T], \mathbb{R}^d) \, \Bigl \vert \,\varphi(0) = x_0\, \, S_{0T}^{(\mathcal{K}_i, \mathcal{K}_{\neg i}^{\ast})}(\varphi(t)) < \alpha \Bigl\}. \label{Eq17}
\end{align}
\end{enumerate} 
where $(\mathcal{K}_i, \mathcal{K}_{\neg i}^{\ast}) \in \mathscr{K}$ for each $i=1, 2, \ldots, n$.
\end{lemma}

\begin{lemma} \label{L2}
Let $D_{+\delta}$ denote a $\delta$-neighborhood of $D$ and let $D_{-\delta}$ denote the set of points in $D$ at a distance greater than $\delta$ from the boundary $\partial D$. Then, for sufficiently small $\delta > 0$, the following estimates 
\begin{align}
\inf_{\substack{\varphi(t) \in C_{0T}([0,T], \mathbb{R}^d) \\ \varphi(0)=x_0}} \biggl \{ S_{0T}^{(\mathcal{K}_i, \mathcal{K}_{\neg i}^{\ast})}(\varphi(t)) \, \Bigl \vert \,\varphi(t) \in D_{+\delta} \cup \partial D_{+\delta}, \forall t \in [0, T] \biggr \}, \label{Eq18}
\end{align}
with $(\mathcal{K}_i, \mathcal{K}_{\neg i}^{\ast}) \in \mathscr{K}$ for each $i=1,2, \dots, n$, and 
\begin{align}
\inf_{\substack{\varphi(t) \in C_{0T}([0,T], \mathbb{R}^d)\\ \varphi(0)=x_0}} \biggl \{ S_{0T}^{(\mathcal{K}_i, \mathcal{K}_{\neg i}^{\ast})}(\varphi(t))\, \Bigl \vert \,\varphi(t) \in D_{-\delta} \cup \partial D_{-\delta}, \forall t \in [0, T] \biggr\}, \label{Eq19}
\end{align}
can be made arbitrarily close to each other. Furthermore, the same holds for
\begin{align}
\inf_{\substack{\varphi(t) \in C_{0T}([0,T], \mathbb{R}^d)\\ \varphi(0)=x,\,\varphi(T)=y}} \biggl \{ S_{0T}^{(\mathcal{K}_i, \mathcal{K}_{\neg i}^{\ast})}(\varphi(t))\, \Bigl \vert \,\varphi(t) \in D_{\pm\delta} \cup \partial D_{\pm\delta}, \forall t \in [0, T] \biggr\}, \label{Eq20}
\end{align}
uniformly for any $x, y \in D_{-\delta}$.
\end{lemma}

Before concluding this section, it is worth mentioning that the asymptotic behavior for singularly perturbed eigenvalue is related to the problem of estimating the minimum asymptotic exit rate with which the state trajectories $x^{\epsilon}(t)$ exit from the domain $D$. For example, for sufficiently small $\epsilon > 0$, the asymptotic behavior of the principal eigenvalue for the infinitesimal generator (corresponding to a zero boundary condition on $\partial D$) has been studied in the past (e.g., see \cite{Day83} or \cite{DevFr78} in the context of an asymptotic behavior for the principal eigenfunction; and see \cite{Kif80}, \cite{Kif81} or \cite{Day87} in the context of an asymptotic behavior for the equilibrium density). Specifically, the authors in \cite{Day87} and \cite{Kif80} have provided some results about the regularity properties of the {\em action functional} in connection with the asymptotic behavior of the equilibrium density, where the latter (i.e., the asymptotic behavior of the equilibrium density) is linked with the exit problem from the domain of attraction with an exponentially stable critical point for the stochastically perturbed dynamical system (see also \cite{She91} and \cite{DayDa84}).

\section{Main Results} \label{S2}

In this section, we present our main results -- where we provide a sufficient condition for the existence of a Nash equilibrium point for the noncooperative $n$-player principal eigenvalue game. Specifically, in the course of such a game, each player generates automatically an admissible control strategy in response to the action of other players via the system state trajectory $x^{\epsilon}(t)$ for $t \ge 0$. For example, the $i$th-player can generate an admissible control strategy $u_i(t)=\bigl(\mathcal{K}_i x^{\epsilon}\bigr)(t)$ in response to the actions of other players $u_j^{\ast}(t)= \bigl(\mathcal{K}_j^{\ast} x^{\epsilon}\bigr)(t)$, for $j \neq i$, with $i = 1,2, \ldots, n$ (where the admissible control strategies $\bigl(\widehat{u_i(\cdot), u_{\neg i}^{\ast}(\cdot)}\bigr) \in \prod\nolimits_{i=1}^n \mathcal{U}_i$ {\em for almost all} $t \ge 0$).\footnote{$\bigl(\widehat{u_i(t), u_{\neg i}^{\ast}(t)}\bigr)\triangleq \bigl(u_1^{\ast}(t),\ldots u_{i-1}^{\ast}(t), u_{i}(t), u_{i+1}^{\ast}(t), \ldots u_{n}^{\ast}(t)\bigr)$.} Moreover, for such a game to have a stable Nash equilibrium point (which is also robust to small perturbations in the strategies played by the other players), then each player is required to respond optimally (in some sense of best-response correspondences) to the actions of the other players.

To this end, it will be useful to consider the following criterion functions (cf. Equation~\eqref{Eq12} or Equation~\eqref{Eq10})
\begin{align}
\mathscr{K}_i[\mathcal{X}, \mathcal{U}_i] \ni \mathcal{K}_i \mapsto r_{i}(\mathcal{K}_i, \mathcal{K}_{\neg i}^{\ast}) \in \mathbb{R}_{-} \cup \{-\infty\}, \,\, i =1,2, \ldots, n, \label{Eq21}
\end{align}
with respect to the admissible controls $u_i(\cdot) \in \mathcal{U}_i$ of the form $u_i(t)=\bigl(\mathcal{K}_i x^{\epsilon}\bigr)(t)$, for $t \ge 0$, where $\mathcal{K}_i  \in \mathscr{K}_i[\mathcal{X}, \mathcal{U}_i]$ for each $i =1, 2, \ldots, n$; while the others $\mathcal{K}_{\neg i}^{\ast}$ remain fixed.\footnote{Notice that such criterion function is upper-semicontinuous and always achieves its extremum over a nonempty closed subset $\prod\nolimits_{i=1}^n \mathscr{K}_i[\mathcal{X}, \mathcal{U}_i]$.} On the other hand, under the game-theoretic setting, if we further assume that the $n$-tuple of linear feedback operators $\mathcal{K}^{\ast} \triangleq \bigl(\mathcal{K}_1^{\ast},\mathcal{K}_2^{\ast}, \ldots, \mathcal{K}_n^{\ast}\bigr) \in \mathscr{K}$ is a Nash equilibrium point. Then, when all players play simultaneously their Nash strategies $u_i^{\ast}(t) = \bigl(\mathcal{K}_i^{\ast} x^{\epsilon}\bigr)(t)$ for all $t \ge 0$, for $i=1,2, \ldots, n$, there exist a unique maximum closed invariant set $\Lambda_D^{\mathcal{K}^{\ast}} \subset D \cup \partial D$ that contains all closed invariant sets from the set $D \cup \partial D$ (under the action of the $n$-tuple of equilibrium linear feedback operators with respect to the unperturbed multi-channel dynamical system) i.e., $\Lambda_D^{(\mathcal{K}_i, \mathcal{K}_{\neg i}^{\ast})} \subseteq \Lambda_D^{\mathcal{K}^{\ast}}$, $\forall \mathcal{K}_i \in \mathscr{K}_{i}[\mathcal{X}, \mathcal{U}_i]$, $\forall i \in \{1,2, \ldots, n\}$. Our interest in this section is to investigate the connection between the Nash equilibrium point and the maximum closed invariant set, and provide a sufficient condition for the existence of Nash equilibrium point for the noncooperative $n$-player principal eigenvalue game (where such a game-theoretic setting further involves some information on the infinitesimal generator of the stochastically perturbed multi-channel dynamical system).\footnote{Note that such a connection is also implicitly related with the problem of maximizing the mean exit time for the controlled state trajectory from the region $D$.}

Therefore, more formally, we have the following definition for the Nash equilibrium point (i.e., the $n$-tuple of equilibrium linear feedback operators). 
\begin{definition} \label{DFN2}
The $n$-tuple $(\mathcal{K}_1^{\ast}, \mathcal{K}_2^{\ast}, \ldots, \mathcal{K}_n^{\ast}) \in \mathscr{K}$ (i.e., the $n$-tuple of equilibrium linear feedback operators) is called a Nash equilibrium point for the principal eigenvalue game if it satisfies
\begin{align}
r_i\bigl(\mathcal{K}_i^{\ast}, \mathcal{K}_{\neg i}^{\ast}\bigr) \le r_i\bigl(\mathcal{K}_i, \mathcal{K}_{\neg i}^{\ast}\bigr), \,\, \forall \mathcal{K}_i \in \mathscr{K}_i[\mathcal{X}, \mathcal{U}_i], \, \, \forall i \in \{1, 2, \ldots, n\}, \label{Eq22}
\end{align}
with $\bigl(\mathcal{K}_i, \mathcal{K}_{\neg i}^{\ast}\bigr) \in \mathscr{K}$ for each $i=1,2, \ldots n$.
\end{definition}

\begin{proposition} \label{P1}
If, for some $x_0 \in D$, one of the following conditions hold
\begin{enumerate} [(i)]
\item
\begin{align}
\lim_{\epsilon \rightarrow 0} \, \limsup_{T \rightarrow \infty} \frac{\epsilon} {T} \log \mathbb{P}_{\epsilon, i}^{(\mathcal{K}_i, \mathcal{K}_{\neg i}^{\ast})} \bigl\{\tau_{D, i}^{\epsilon} > T \bigr\} > -\infty, \label{Eq23}
\end{align}
\item
\begin{align}
\lim_{\epsilon \rightarrow 0} \,\mathbb{E}_{\epsilon, i}^{(\mathcal{K}_i, \mathcal{K}_{\neg i}^{\ast})} \bigl\{\tau_{D, i}^{\epsilon}\bigr\} = \infty, \label{Eq24}
\end{align}
\end{enumerate}
where $\tau_{D, i}^{\epsilon}$ is the exit-time with respect to the $i$th-channel and $(\mathcal{K}_i, \mathcal{K}_{\neg i}^{\ast}) \in \mathscr{K}$ for each $i=1,2, \dots, n$. Then, the maximum closed invariant set $\Lambda_D^{(\mathcal{K}_i, \mathcal{K}_{\neg i}^{\ast})}$ for the dynamical system in Equation~\eqref{Eq5} is nonempty.
\end{proposition}

Then, we have the following proposition which is a direct consequence of Proposition~\ref{P1}.
\begin{proposition} \label{P2}
Suppose that the class of linear feedback operators $\mathscr{K}$ is nonempty. Then, the principal eigenvalue $\lambda_{\epsilon, i}^{(\mathcal{K}_i, \mathcal{K}_{\neg i}^{\ast})}$ corresponding to the infinitesimal generator $\mathcal{L}_{\epsilon}^{(\mathcal{K}_i, \mathcal{K}_{\neg i}^{\ast})}$ with zero boundary condition on $\partial D$ satisfies 
\begin{align}
\lambda_{\epsilon, i}^{(\mathcal{K}_i, \mathcal{K}_{\neg i}^{\ast})} = \epsilon^{-1} r_{i}(\mathcal{K}_i, \mathcal{K}_{\neg i}^{\ast}) + o(\epsilon^{-1}) \quad \text{as} \quad \epsilon \rightarrow 0, \label{Eq25}
\end{align}
where
\begin{align}
r_{i}(\mathcal{K}_i, \mathcal{K}_{\neg i}^{\ast}) = \limsup_{T \rightarrow \infty} \,\inf_{\substack{\varphi(t) \in C_{0T}([0,T], \mathbb{R}^d)\\ \varphi(0)=x_0}} \frac{1}{T} \biggl \{ S_{0T}^{(\mathcal{K}_i, \mathcal{K}_{\neg i}^{\ast})}(\varphi(t)) \, \Bigl \vert \,\varphi(t) \in D \cup \partial D, \forall t \in [0, T] \biggr\}, \label{Eq26}
\end{align}
with
\begin{align*}
 S_{0T}^{(\mathcal{K}_i, \mathcal{K}_{\neg i}^{\ast})}(\varphi(t)) =  \frac{1}{2} \int_{0}^{T}\biggl\Vert\frac{d\varphi(t)}{dt} - \bigl(A + \bigl(B, (\mathcal{K}_i, \mathcal{K}_{\neg i}^{\ast})\bigr)\varphi\bigr)(t)\biggr\Vert^2dt
 \end{align*}
and $(\mathcal{K}_i, \mathcal{K}_{\neg i}^{\ast}) \in \mathscr{K}$ for each $i = 1, 2, \ldots, n$.
\end{proposition}

The following proposition provides a sufficient condition for the existence of a Nash equilibrium point for the noncooperative $n$-player principal eigenvalue game.
\begin{proposition} \label{P3}
Suppose that the mapping $\mathscr{K}_i[\mathcal{X}, \mathcal{U}_i] \ni \mathcal{K}_i \mapsto r_{i}(\mathcal{K}_i, \mathcal{K}_{\neg i}^{\ast}) \in \mathbb{R}_{+} \cup \{\infty\}$ is upper-semicontinuous. Then, there exists at least one Nash equilibrium point that satisfies
\begin{align}
r_{i}(\mathcal{K}_i^{\ast}, \mathcal{K}_{\neg i}^{\ast}) \le r_{i}(\mathcal{K}_i, \mathcal{K}_{\neg i}^{\ast}),\,\,  \forall \mathcal{K}_i \in \mathscr{K}_i[\mathcal{X}, \mathcal{U}_i], \,\,  \forall i \in \{1, 2, \dots, n\}, \label{Eq27}
\end{align}
where
\begin{align}
\mathcal{K}_{i}^{\ast} \in \argmin_{\mathcal{K}_i \in \mathscr{K}_{i}[\mathcal{X}, \mathcal{U}_i]} \Biggr\{ \limsup_{T \rightarrow \infty} \,\inf_{\substack{\varphi(t) \in C_{0T}([0,T], \mathbb{R}^d)\\ \varphi(0)=x_0}} \frac{1}{T} \biggl \{ S_{0T}^{(\mathcal{K}_i, \mathcal{K}_{\neg i}^{\ast})}(\varphi(t))& \, \Bigl \vert \, \varphi(t) \in D \cup \partial D,\, \notag \\
& \forall t \in [0, T] \biggr\} \Biggl\}, \label{Eq28}
\end{align}
with $(\mathcal{K}_i, \mathcal{K}_{\neg i}^{\ast})\in \mathscr{K}$ for each $i=1,2, \ldots n$.

Furthermore, the maximum closed invariant set $\Lambda_D^{\mathcal{K}^{\ast}} \subset D \cup \partial D$ (under the action of the $n$-tuple of equilibrium linear feedback operators $(\mathcal{K}_1^{\ast}, \mathcal{K}_2^{\ast}, \ldots, \mathcal{K}_n^{\ast}) \in \mathscr{K}$ with respect to the unperturbed multi-channel dynamical system) satisfies
\begin{align}
\Lambda_D^{(\mathcal{K}_i, \mathcal{K}_{\neg i}^{\ast})} \subset \Lambda_D^{\mathcal{K}^{\ast}},\,\, \forall \mathcal{K}_i \in \mathscr{K}_{i}[\mathcal{X}, \mathcal{U}_i],\,\, \forall i \in \{1, 2, \dots, n\}, \label{Eq29}
\end{align}
with $(\mathcal{K}_i, \mathcal{K}_{\neg i}^{\ast})\in \mathscr{K}$ for each $i=1,2, \ldots n$.
\end{proposition}

\section{Proof of the Main Results} \label{S3}
\subsection{Proof of Proposition~\ref{P1}} \label{S3(1)}
For a fixed $i \in \{1, 2, \ldots, n\}$, suppose that the maximum closed invariant set $\Lambda_D^{(\mathcal{K}_i, \mathcal{K}_{\neg i}^{\ast})}$, with $(\mathcal{K}_i, \mathcal{K}_{\neg i}^{\ast}) \in \mathscr{K}$, is empty. Then, there exists an open bounded domain $\tilde{D} \supset D \cup \partial D$ such that the corresponding set $\Lambda_{\tilde{D}}^{(\mathcal{K}_i, \mathcal{K}_{\neg i}^{\ast})}$  is also empty.

Note that it is easy to check that if $D_{2} \subset D_{1}$, then $\Lambda_{D_{2}}^{(\mathcal{K}_i, \mathcal{K}_{\neg i}^{\ast})} \subset \Lambda_{D_{1}}^{(\mathcal{K}_i, \mathcal{K}_{\neg i}^{\ast})}$. Take the following sequence $\bigl\{D_{m}\bigr\}$ of open domains such that
\begin{align}
D_{1} \supset D_{2} \supset D_{3} \supset \cdots \quad \text{and} \quad \bigcap\nolimits_{m \ge 1} D_{m} = D \cup \partial D. \label{Eq30}
\end{align}
If $\Lambda_{D_{m}}^{(\mathcal{K}_i, \mathcal{K}_{\neg i}^{\ast})} \neq \varnothing$ for all $m \ge 1$, then
\begin{align}
\Lambda = \bigcap\nolimits_{m \ge 1} \Lambda_{D_{m}}^{(\mathcal{K}_i, \mathcal{K}_{\neg i}^{\ast})}. \label{Eq31}
\end{align}
Moreover, since $\Lambda_{D_{m}}^{(\mathcal{K}_i, \mathcal{K}_{\neg i}^{\ast})}$ is closed, we have
\begin{align}
\Lambda_{D_{1}}^{(\mathcal{K}_i, \mathcal{K}_{\neg i}^{\ast})} \supset \Lambda_{D_{2}}^{(\mathcal{K}_i, \mathcal{K}_{\neg i}^{\ast})} \supset \Lambda_{D_{3}}^{(\mathcal{K}_i, \mathcal{K}_{\neg i}^{\ast})} \supset  \cdots. \label{Eq32}
\end{align}
Note that  $\Lambda$ is an invariant closed set with respect to the unperturbed multi-channel dynamical system and $\Lambda \supset D \cup \partial D$. Thus, $\varnothing \neq \Lambda \subset \Lambda_D^{(\mathcal{K}_i, \mathcal{K}_{\neg i}^{\ast})}$. This contradicts our earlier assumption. Then, for some $m_0 \ge 1$, we have 
\begin{align}
\Lambda_{D_{m_0}}^{(\mathcal{K}_i, \mathcal{K}_{\neg i}^{\ast})} = \varnothing. \label{Eq33}
\end{align}
Let $\tilde{D} = D_{m_0}$ and, for any $x_0 \in \tilde{D} \cup \partial \tilde{D}$, let us introduce the following
\begin{align}
\tau_{\tilde{D}, i}^{0} = \inf \bigl\{ t > 0 \, \bigl\vert \, x^{0}(t) \notin  \tilde{D} \cup \partial \tilde{D} \bigr\}, \label{Eq34}
\end{align}
with respect to $(\mathcal{K}_i, \mathcal{K}_{\neg i}^{\ast}) \in \mathscr{K}$ for each $i \in \{1, 2, \ldots, n\}$. Then, we can show that
\begin{align}
\tau_{\tilde{D}, i}^{0} < \infty, \label{Eq35}
\end{align}
for any $x_0 \in \tilde{D} \cup \partial \tilde{D}$. Note that, if $\tau_{\tilde{D}, i}^{0} =\infty$, then $\phi\bigl(t; 0, x_0, (\widehat{\mathcal{K}_i x^{0}, \mathcal{K}_{\neg i}^{\ast}x^{0}}) (t)\bigr) \in \tilde{D} \cup \partial \tilde{D}$ for all $t \ge 0$. Then, for some sequence $t_m \rightarrow \infty$ and a point $y \in \tilde{D} \cup \partial \tilde{D}$, we have
\begin{align}
\phi\bigl(t_m; 0, x_0, (\widehat{\mathcal{K}_i x^{0}, \mathcal{K}_{\neg i}^{\ast}x^{0}})(t_n)\bigr) \rightarrow y \quad \text{as} \quad t_m \rightarrow \infty \label{Eq36}
\end{align}
and
\begin{align}
\phi\bigl((t_m+t); 0, x_0, (\widehat{\mathcal{K}_i x^{0}, \mathcal{K}_{\neg i}^{\ast}x^{0}})(t_m+t)\bigr) \rightarrow \phi\bigl((t; t_m, y, (\widehat{\mathcal{K}_i x^{0}, \mathcal{K}_{\neg i}^{\ast}x^{0}})(t)\bigr), \label{Eq37}
\end{align}
for any $t \in [0, \infty)$.

Thus, if $ \phi\bigl(t; t_m, y, (\widehat{\mathcal{K}_i x^{0}, \mathcal{K}_{\neg i}^{\ast}x^{0}})(t)\bigr) \in \tilde{D} \cup \partial \tilde{D}$ for all $t \in [0, \infty)$, then we have the following
\begin{align}
\Bigl\{\phi\bigl((t; t_m, y, (\widehat{\mathcal{K}_i x^{0}, \mathcal{K}_{\neg i}^{\ast}x^{0}})(t)\bigr), \,\, t \ge 0 \Bigr\} \subset \Lambda_{\tilde{D}}^{(\mathcal{K}_i, \mathcal{K}_{\neg i}^{\ast})} = \varnothing, \label{Eq38}
\end{align}
which show that $\tau_{\tilde{D}, i}^{0}$ is finite.

Note that, from upper-semicontinuity of $\tau_{\tilde{D}, i}^{0}$, we have
\begin{align}
\tilde{T} = \sup_{x_0 \in \tilde{D} \cup \partial \tilde{D}} \tau_{\tilde{D}, i}^{0} < \infty. \label{Eq39}
\end{align}
Moreover, for any $\delta > 0$, let\footnote{Here the diffusion process $x^{\epsilon}(t)$ is described by the following stochastic differential equation
\begin{align*}
d x^{\epsilon}(t) = A x^{\epsilon}(t) dt + B_i (\mathcal{K}_i x^{\epsilon})(t) dt +  \sum\nolimits_{j \neq i} B_j (\mathcal{K}_j^{\ast}x^{\epsilon})(t) + \sqrt{\epsilon} \sigma(x^{\epsilon}(t))dW(t), \,\, x^{\epsilon}(0) = x_0.
\end{align*}}
\begin{align}
\lim_{\epsilon \rightarrow 0} \sup_{x_0 \in D} \mathbb{P}_{\epsilon, i}^{(\mathcal{K}_i, \mathcal{K}_{\neg i}^{\ast})} \bigl\{ \rho(x^{\epsilon}(t), \phi\bigl(t; 0, x_0, (\widehat{\mathcal{K}_i x^{0}, \mathcal{K}_{\neg i}^{\ast}x^{0}})(t)\bigr)) > \delta \bigr\} = 0, \,\, t \ge 0. \label{Eq40}
\end{align}
From Equations~\eqref{Eq34}--\eqref{Eq41}, we have
\begin{align}
 \mathbb{P}_{\epsilon, i}^{(\mathcal{K}_i, \mathcal{K}_{\neg i}^{\ast})} \bigl\{ \tau_{D, i}^{\epsilon} > \tilde{T} \bigr\} \rightarrow 0 \quad \text{as} \quad \epsilon \rightarrow 0. \label{Eq41}
\end{align}
Then, using the Markov property, we have
\begin{align}
 \mathbb{P}_{\epsilon, i}^{(\mathcal{K}_i, \mathcal{K}_{\neg i}^{\ast})} \bigl\{ \tau_{D, i}^{\epsilon} > \ell \tilde{T} \bigr\} &= \mathbb{E}_{\epsilon, i}^{(\mathcal{K}_i, \mathcal{K}_{\neg i}^{\ast})} \chi_{\tau_{D, i}^{\epsilon} > \tilde{T}} \mathbb{E}_{\epsilon, i}^{(\mathcal{K}_i, \mathcal{K}_{\neg i}^{\ast})} \chi_{\tau_{D, i}^{\epsilon} > \tilde{T}} \cdots \mathbb{E}_{\epsilon, i}^{(\mathcal{K}_i, \mathcal{K}_{\neg i}^{\ast})} \chi_{\tau_{D, i}^{\epsilon} > \tilde{T}}, \notag \\
 & \le \Bigl(\mathbb{P}_{\epsilon, i}^{(\mathcal{K}_i, \mathcal{K}_{\neg i}^{\ast})} \bigl\{ \tau_{D, i}^{\epsilon} > \tilde{T} \bigr\} \Bigr)^{\ell}, \label{Eq42}
\end{align}
where $\chi_{\mathcal{A}}$ is the indicator for the event $\mathcal{A}$.

Since $\mathbb{P}_{\epsilon, i}^{(\mathcal{K}_i, \mathcal{K}_{\neg i}^{\ast})} \bigl\{ \tau_{D, i}^{\epsilon} > T \bigr\}$ decreases in $T$, then we have
\begin{align}
\limsup_{T \rightarrow \infty} \frac{1}{T} \log \mathbb{P}_{\epsilon, i}^{(\mathcal{K}_i, \mathcal{K}_{\neg i}^{\ast})} \bigl\{ \tau_{D, i}^{\epsilon} > T \bigr\} \le  \frac{1}{\tilde{T}} \log \mathbb{P}_{\epsilon, i}^{(\mathcal{K}_i, \mathcal{K}_{\neg i}^{\ast})} \bigl\{ \tau_{D, i}^{\epsilon} > \tilde{T} \bigr\}. \label{Eq43}
\end{align}
Taking into account Equation~\eqref{Eq42}, then, for any $x_0 \in D$, we have the following
\begin{align}
 \limsup_{T \rightarrow \infty} \frac{1}{T} \log \mathbb{P}_{\epsilon, i}^{(\mathcal{K}_i, \mathcal{K}_{\neg i}^{\ast})} \bigl\{\tau_{D, i}^{\epsilon} > T \bigr\}  \rightarrow -\infty \quad \text{as} \quad \epsilon \rightarrow 0. \label{Eq44}
\end{align}
Hence, our assumption that $\Lambda_D^{(\mathcal{K}_i, \mathcal{K}_{\neg i}^{\ast})} = \varnothing$ is inconsistent.

To proof the part (ii), notice that
\begin{align}
 \mathbb{E}_{\epsilon, i}^{(\mathcal{K}_i, \mathcal{K}_{\neg i}^{\ast})} \bigl\{ \tau_{D, i}^{\epsilon} \bigr\} \le \tilde{T}\sum\nolimits_{\ell=1}^{\infty} \mathbb{P}_{\epsilon, i}^{(\mathcal{K}_i, \mathcal{K}_{\neg i}^{\ast})} \bigl\{ \tau_{D, i}^{\epsilon} > (\ell-1)\tilde{T} \bigr\}. \label{Eq45}
\end{align}
Assumption  $\Lambda_D^{(\mathcal{K}_i, \mathcal{K}_{\neg i}^{\ast})} = \varnothing$ gives, in view of Equations~\eqref{Eq41}, \eqref{Eq42} and \eqref{Eq45} for sufficiently small $\epsilon > 0$ and for any $x_0 \in D$, that
\begin{align}
 \mathbb{E}_{\epsilon, i}^{(\mathcal{K}_i, \mathcal{K}_{\neg i}^{\ast})} \bigl\{ \tau_{D, i}^{\epsilon} \bigr\} &\le \tilde{T}\sum\nolimits_{\ell=1}^{\infty} \mathbb{P}_{\epsilon, i}^{(\mathcal{K}_i, \mathcal{K}_{\neg i}^{\ast})} \bigl\{ \tau_{D, i}^{\epsilon} > (m-1)\tilde{T} \bigr\}, \notag \\
 &< \infty, \label{Eq46}
\end{align}
which contradicts with Equation~\eqref{Eq23}. This completes the proof of Proposition~\ref{P1}. \qed

\subsection{Proof of Proposition~\ref{P2}} \label{S3(2)}
For a fixed $i \in \{1, 2, \ldots, n\}$, suppose that $r_{i}(\mathcal{K}_i, \mathcal{K}_{\neg i}^{\ast})$, with $(\mathcal{K}_i, \mathcal{K}_{\neg i}^{\ast}) \in \mathscr{K}$, exists.\footnote{Note that the existence of such a limit for $r_{i}(\mathcal{K}_i, \mathcal{K}_{\neg i}^{\ast})$ can be easily established (e.g., see \cite{VenFre70}).} Then, using Lemma~\ref{L2}, one can show that $r_{i}(\mathcal{K}_i, \mathcal{K}_{\neg i}^{\ast})$ also satisfies the following
\begin{align}
r_{i}(\mathcal{K}_i, \mathcal{K}_{\neg i}^{\ast}) =  \sup_{x, y \in D} \Biggl\{ \limsup_{T \rightarrow \infty} \, \inf_{\substack{\varphi(t) \in C_{0T}([0,T], \mathbb{R}^d)\\ \varphi(0)=x,\,\varphi(T)=y}} \frac{1}{T} \biggl \{ &S_{0T}^{(\mathcal{K}_i, \mathcal{K}_{\neg i}^{\ast})}(\varphi(t))\, \Bigl \vert \,\varphi(t) \in D \cup \partial D, \notag \\
 &\quad \quad\quad\quad \forall t \in [0, T] \biggr\} \Biggr\}. \label{Eq47}
\end{align}
Next, let us show that, for sufficiently small $\epsilon > 0$, $\mathbb{E}_{\epsilon, i}^{(\mathcal{K}_i, \mathcal{K}_{\neg i}^{\ast})}\bigl\{\exp(\epsilon^{-1} R \tau_{D, i}^{\epsilon})\bigr\}$ tends to infinity, when $R > r_{i}(\mathcal{K}_i, \mathcal{K}_{\neg i}^{\ast})$. If we choose a positive $\varkappa$ which is smaller than $(R - r_{i}(\mathcal{K}_i, \mathcal{K}_{\neg i}^{\ast}))/3$ so that
\begin{align}
 \sup_{x, y \in D} \, \inf_{\substack{\varphi(t) \in C_{0T}([0,T], \mathbb{R}^d)\\ \varphi(0)=x,\,\varphi(T)=y}} \frac{1}{T} \biggl \{ S_{0T}^{(\mathcal{K}_i, \mathcal{K}_{\neg i}^{\ast})}(\varphi(t))\, \Bigl \vert \,\varphi(t) \in D \cup \partial D, \forall t \in [0, T] \biggr\} \notag \\
 < r_{i}(\mathcal{K}_i, \mathcal{K}_{\neg i}^{\ast}) + \varkappa, \label{Eq48}
\end{align}
and, for sufficiently small $\delta > 0$,
\begin{align}
 \inf_{\substack{\varphi(t) \in C_{0T}([0,T], \mathbb{R}^d)\\ \varphi(0)=x,\,\varphi(T)=y}} \biggl \{ S_{0T}^{(\mathcal{K}_i, \mathcal{K}_{\neg i}^{\ast})}(\varphi(t))\, \Bigl \vert \,\varphi(t) \in D_{-\delta} \cup \partial D_{-\delta}, \forall t \in [0, T] \biggr\} \notag \\
 < T\bigl(r_{i}(\mathcal{K}_i, \mathcal{K}_{\neg i}^{\ast}) + 2\varkappa\bigr), \label{Eq49}
\end{align}
for all $x, y \in D_{-\delta}$. Then, if we further let $\alpha = T\bigl(r_{i}(\mathcal{K}_i, \mathcal{K}_{\neg i}^{\ast}) + 2 \varkappa\bigr)$ and $\gamma = T\varkappa$, from Lemma~\ref{L1}, there exits an $\epsilon_0 > 0$ such that
\begin{align}
 S_{0T}^{(\mathcal{K}_i, \mathcal{K}_{\neg i}^{\ast})}(\varphi(t)) \le T\bigl(r_{i}(\mathcal{K}_i, \mathcal{K}_{\neg i}^{\ast}) + 2 \varkappa\bigr), \quad \varphi(t) \in D_{-\delta} \cup \partial D_{-\delta}, \forall t \in [0, T], \label{Eq50}
\end{align}
for any $x, y \in D_{-\delta}$; and, moreover, we have the following probability estimate
\begin{align}
\mathbb{P}_{\epsilon, i}^{(\mathcal{K}_i, \mathcal{K}_{\neg i}^{\ast})} \Bigl\{ \tau_{D_{-\delta}, i}^{\epsilon} > T,\, x^{\epsilon}(\tau_{D_{-\delta}, i}^{\epsilon}) \in D_{-\delta} \Bigr\} & \ge \mathbb{P}_{\epsilon, i}^{(\mathcal{K}_i, \mathcal{K}_{\neg i}^{\ast})} \Bigl\{\rho_{0T}\bigl(x^{\epsilon}(t), \varphi(t)\bigr) < \delta \Bigr\}, \notag \\ 
&\ge \exp \bigl(-\epsilon^{-1}\bigl(S_{0T}^{(\mathcal{K}_i, \mathcal{K}_{\neg i}^{\ast})}(\varphi(t)) + \gamma \bigr)\bigl), \notag \\
&\ge \exp \bigl( -\epsilon^{-1} T \bigl(r_{i}(\mathcal{K}_i, \mathcal{K}_{\neg i}^{\ast}) + 3\varkappa \bigr)\bigl),\notag \\
& \quad\quad\quad\quad\quad\quad\quad \forall  \epsilon \in (0, \epsilon_0), \label{Eq51}
\end{align}
where $\varphi(T) \in D_{-2\delta}$.

Let us define the following random events 
\begin{align}
\mathcal{A}_{\ell} = \Bigl\{ \tau_{D_{-\delta}, i}^{\epsilon} > \ell T,\, x^{\epsilon}(\ell T) \in D_{-\delta}\Bigl\}, \label{Eq52}
\end{align}
for $\ell \in \mathbb{N}_{+} \cup \{0\}$. Then, from the Markov property, we have
\begin{align}
\mathbb{P}_{\epsilon, i}^{(\mathcal{K}_i, \mathcal{K}_{\neg i}^{\ast})} \bigl\{ \mathcal{A}_{\ell} \bigl\} & \ge \mathbb{E}_{\epsilon, i}^{(\mathcal{K}_i, \mathcal{K}_{\neg i}^{\ast})}  \chi_{\mathcal{A}_{\ell-1}} \mathbb{P}_{(\epsilon, x_{(\ell-1)T})}^{(\mathcal{K}_i, \mathcal{K}_{\neg i}^{\ast})}  \bigl\{ \mathcal{A}_1 \bigr\},\notag  \\
& \ge \mathbb{P}_{\epsilon, i}^{(\mathcal{K}_i, \mathcal{K}_{\neg i}^{\ast})} \bigl\{ \mathcal{A}_{\ell-1} \bigr\}  \inf_{y \in D_{-\delta}} \mathbb{P}_{(\epsilon, y)}^{(\mathcal{K}_i, \mathcal{K}_{\neg i}^{\ast})} \bigl\{ \mathcal{A}_1 \bigr\},\notag  \\
& \ge \exp \bigl( - \epsilon^{-1} \ell T\bigl(r_{i}(\mathcal{K}_i, \mathcal{K}_{\neg i}^{\ast}) + 3\varkappa \bigr)\bigr) \quad \forall  \epsilon \in (0, \epsilon_0). \label{Eq53}
\end{align}
Note that, for an arbitrary $\ell$, we have the following
\begin{align}
\mathbb{E}_{\epsilon, i}^{(\mathcal{K}_i, \mathcal{K}_{\neg i}^{\ast})} \bigl\{ \exp \bigl( \epsilon^{-1} R \tau_{D_{-\delta}, i}^{\epsilon} \bigl) \bigl\} & \ge \exp \bigl( \epsilon^{-1} R \ell T \bigl) \mathbb{P}_{\epsilon, i}^{(\mathcal{K}_i, \mathcal{K}_{\neg i}^{\ast})} \bigl\{ \tau_{D_{-\delta}, i}^{\epsilon} >  \ell T \bigr\},\notag  \\
& \ge \exp \bigl( - \epsilon^{-1} \ell T\bigl(R - r_{i}(\mathcal{K}_i, \mathcal{K}_{\neg i}^{\ast}) - 3\varkappa \bigr)\bigr), \quad \forall  \epsilon \in (0, \epsilon_0), \label{Eq54}
\end{align}
which tends to infinity as $\ell \rightarrow \infty$, i.e., $\mathbb{E}_{\epsilon, i}^{(\mathcal{K}_i, \mathcal{K}_{\neg i}^{\ast})} \bigl\{ \exp \bigl(\epsilon^{-1} R \tau_{D_{-\delta}, i}^{\epsilon} \bigl) \bigl\}=\infty$.

On the other hand, let us show that if $R < r_{i}(\mathcal{K}_i, \mathcal{K}_{\neg i}^{\ast})$, then, for sufficiently small $\epsilon > 0$, $\mathbb{E}_{\epsilon, i}^{(\mathcal{K}_i, \mathcal{K}_{\neg i}^{\ast})} \bigl\{ \exp \bigl(\epsilon^{-1} R \tau_{D_{-\delta}, i}^{\epsilon} \bigl) \bigl\} < \infty$. For $\varkappa < (R - r_{i}(\mathcal{K}_i, \mathcal{K}_{\neg i}^{\ast}))/3$, let us choose $\delta$ so that
\begin{align}
\inf_{\substack{\varphi(t) \in C_{0T}([0,T], \mathbb{R}^d) \\ \varphi(0)=x_0}} \biggl \{ S_{0T}^{(\mathcal{K}_i, \mathcal{K}_{\neg i}^{\ast})}(\varphi(t)) \, \Bigl \vert \,\varphi(t) \in D_{+\delta} \cup \partial D_{+\delta}, \forall t \in [0, T] \biggr \} \notag \\
> T\bigl(r_{i}(\mathcal{K}_i, \mathcal{K}_{\neg i}^{\ast}) - 2 \varkappa \bigr). \label{Eq55}
\end{align}
From Lemma~\ref{L1}, with $\alpha = T\bigl(r_{i}(\mathcal{K}_i, \mathcal{K}_{\neg i}^{\ast}) - 2 \varkappa \bigr)$ and $\gamma T\varkappa$, there exists an $\epsilon_0 > 0$ such that the distance between the set of functions $\psi(t)$, for $0 \le t \le T$, entirely lying in $D$ and any of the sets $\Phi_{x_0, \alpha}^{(\mathcal{K}_i, \mathcal{K}_{\neg i}^{\ast})}$ is at least a distance $\delta$; and, hence, we have the following probability estimate
\begin{align}
\mathbb{P}_{\epsilon, i}^{(\mathcal{K}_i, \mathcal{K}_{\neg i}^{\ast})} \Bigl\{ \tau_{D, i}^{\epsilon} > T \Bigr\} & \le \mathbb{P}_{\epsilon, i}^{(\mathcal{K}_i, \mathcal{K}_{\neg i}^{\ast})} \Bigl\{d_{0T}\bigl(x^{\epsilon}(t), \Phi_{x_0, \alpha}^{(\mathcal{K}_i, \mathcal{K}_{\neg i}^{\ast})} \bigr) \ge \delta \Bigr\},\notag \\ 
&\le \exp \bigl(-\epsilon^{-1} T \bigl(r_{i}(\mathcal{K}_i, \mathcal{K}_{\neg i}^{\ast}) - 3\varkappa \bigr)\bigl), \quad \forall  \epsilon \in (0, \epsilon_0), \label{Eq56}
\end{align}
for any $x \in D$.

Then, using the Markov property, we have the following
\begin{align}
\mathbb{P}_{\epsilon, i}^{(\mathcal{K}_i, \mathcal{K}_{\neg i}^{\ast})} \Bigl\{ \tau_{D, i}^{\epsilon} > \ell T \Bigr\} &\le \exp \bigl( -\epsilon^{-1} \ell T \bigl(r_{i}(\mathcal{K}_i, \mathcal{K}_{\neg i}^{\ast}) - 3\varkappa \bigr)\bigl), \quad \forall  \epsilon \in (0, \epsilon_0), \label{Eq57}
\end{align}
and
\begin{align}
&\mathbb{E}_{\epsilon, i}^{(\mathcal{K}_i, \mathcal{K}_{\neg i}^{\ast})} \bigl\{ \exp \bigl( \epsilon^{-1} R \tau_{D, i}^{\epsilon} \bigl) \bigl\} \notag \\
&\quad \quad \quad \quad \le \sum\nolimits_{\ell=0}^{\infty} \exp \bigl( \epsilon^{-1} R (\ell +1) T \bigl) \mathbb{P}_{\epsilon, i}^{(\mathcal{K}_i, \mathcal{K}_{\neg i}^{\ast})} \bigl\{\ell T <  \tau_{D, i}^{\epsilon} \le (\ell +1)T \bigr\},\notag  \\
&\quad \quad \quad \quad \le \sum\nolimits_{\ell=0}^{\infty} \exp \bigl( \epsilon^{-1} R (\ell +1) T \bigl) \mathbb{P}_{\epsilon, i}^{(\mathcal{K}_i, \mathcal{K}_{\neg i}^{\ast})} \bigl\{ \tau_{D, i}^{\epsilon} > \ell T \bigr\},\notag  \\
&\quad \quad \quad \quad \le \sum\nolimits_{\ell=0}^{\infty} \exp \bigl( \epsilon^{-1} R T \bigl) \exp \bigl( -\epsilon^{-1} \ell T\bigl(R - r_{i}(\mathcal{K}_i, \mathcal{K}_{\neg i}^{\ast}) + 3\varkappa \bigr)\bigr), \,\, \forall \epsilon \in (0, \epsilon_0), \label{Eq58}
\end{align}
which converges to a finite value, i.e., $ \mathbb{P}_{\epsilon, i}^{(\mathcal{K}_i, \mathcal{K}_{\neg i}^{\ast})} \bigl\{ \exp \bigl( \epsilon^{-1} R \tau_{D, i}^{\epsilon} \bigl) \bigl\} < \infty$. Hence, $r_{i}(\mathcal{K}_i, \mathcal{K}_{\neg i}^{\ast})$ is a boundary for which $\mathbb{E}_{\epsilon, i}^{(\mathcal{K}_i, \mathcal{K}_{\neg i}^{\ast})} \bigl\{\exp(\epsilon^{-1} r_{i}(\mathcal{K}_i, \mathcal{K}_{\neg i}^{\ast}) \tau_{D, i}^{\epsilon})\bigr\}$ is finite. Then, from Equation~\eqref{Eq35} (cf. Equation\eqref{Eq31}), we have
\begin{align}
-\frac{1}{T} \log \mathbb{P}_{\epsilon, i}^{(\mathcal{K}_i, \mathcal{K}_{\neg i}^{\ast})} \Bigl\{ \tau_{D, i}^{\epsilon} > T \Bigr\} \le \epsilon^{-1}\bigl(r_{i}(\mathcal{K}_i, \mathcal{K}_{\neg i}^{\ast}) - 3\varkappa \bigr), \quad \forall  \epsilon \in (0, \epsilon_0), \label{Eq59}
\end{align}
for any $x \in D$, where the left side tends to the principal eigenvalue $\lambda_{\epsilon, i}^{(\mathcal{K}_i, \mathcal{K}_{\neg i}^{\ast})}$ as $T \rightarrow \infty$. This completes the proof of Proposition~\ref{P2}. \qed

\subsection{Proof of  Proposition~\ref{P3}} \label{S3(3)}
To prove this proposition, we use the Ekeland's variational principle for equilibrium problems (e.g., see \cite{AubEk84}). To this end, for some $x_0 \in D$, let us introduce the following auxiliary mapping $\varrho \colon \prod\nolimits_{i=1}^n \mathscr{K}_i[\mathcal{X}, \mathcal{U}_i] \times \prod\nolimits_{i=1}^n \mathscr{K}_i[\mathcal{X}, \mathcal{U}_i] \to \mathbb{R}_{-}\cup\{-\infty\}$, i.e.,
\begin{align}
\bigl(\mathcal{K}^{\ast}, \mathcal{K} \bigr) \mapsto \varrho(\mathcal{K}^{\ast}, \mathcal{K}) \triangleq \sum\nolimits_{i=1}^n \bigr \{r_{i}(\mathcal{K}_i^{\ast}, \mathcal{K}_{\neg i}^{\ast}) - r_{i}(\mathcal{K}_i, \mathcal{K}_{\neg i}^{\ast}) \bigl \}, \label{Eq60}
\end{align}
which is lower-semicontinuous with respect to $\mathcal{K}=(\mathcal{K}_i, \mathcal{K}_{\neg i}) \in \prod\nolimits_{i=1}^n \mathscr{K}_i[\mathcal{X}, \mathcal{U}_i]$ and it also satisfies the following
\begin{align}
\varrho(\mathcal{K}^{\ast}, \mathcal{K}) \le \varrho(\mathcal{K}^{\ast}, \tilde{\mathcal{K}}) + \varrho(\tilde{\mathcal{K}}, \mathcal{K}), \,\, \forall \mathcal{K}^{\ast}, \mathcal{K}, \tilde{\mathcal{K}} \in \prod\nolimits_{i=1}^n \mathscr{K}_i[\mathcal{X}, \mathcal{U}_i], \label{Eq61}
\end{align}
with $\tilde{\mathcal{K}}=(\tilde{\mathcal{K}}_i, \tilde{\mathcal{K}}_{\neg i}^{\ast}) \in \prod\nolimits_{i=1}^n \mathscr{K}_i[\mathcal{X}, \mathcal{U}_i]$. Moreover, for each $\ell \in \mathbb{N}_{+}$, if $\mathcal{K}^{(\ell)}=(\mathcal{K}_i^{(\ell)}, {\mathcal{K}_{\neg i}^{\ast}}^{(\ell)}) \in \prod\nolimits_{i=1}^n \mathscr{K}_i[\mathcal{X}, \mathcal{U}_i]$ is an $\varepsilon$-equilibrium point.\footnote{Note that if $\Vert \mathcal{K}^{\ast} - \mathcal{K}^{(\ell)}\Vert_{\prod\nolimits_{i=1}^n \mathscr{K}_i[\mathcal{X}, \mathcal{U}_i]}^2 \le \varepsilon$ for sufficiently small $\varepsilon > 0$, then we call $\mathcal{K}^{(\ell)}$ an {\em epsilon-equilibrium}, i.e., a near-Nash equilibrium point (noting that $\mathcal{K}^{\ast}$ is a Nash-equilibrium).} Then, we have
\begin{align}
\varrho(\mathcal{K}^{\ast}, \mathcal{K}^{(\ell)}) \ge - \varepsilon \bigl\Vert \mathcal{K}^{\ast} - \mathcal{K}^{(\ell)} \bigr\Vert_{\prod\nolimits_{i=1}^n \mathscr{K}_i[\mathcal{X}, \mathcal{U}_i]}^2, \,\, \forall \mathcal{K}^{(\ell)} \in \prod\nolimits_{i=1}^n \mathscr{K}_i[\mathcal{X}, \mathcal{U}_i]. \label{Eq62}
\end{align}
Notice that $\varrho(\mathcal{K}^{\ast},\,.)$ is upper-semicontinuous for every $\mathcal{K}$ from the closed set $\prod\nolimits_{i=1}^n \mathscr{K}_i[\mathcal{X}, \mathcal{U}_i]$, then we can choose a subsequence $\bigl\{\mathcal{K}^{(\ell_k)}\bigr\}$ of $\bigl\{\mathcal{K}^{(\ell)}\bigr\}$ such that $\mathcal{K}^{(\ell_k)} \rightarrow \mathcal{K}^{\ast}$ as $k \rightarrow \infty$. Hence, we have
\begin{align}
\liminf_{\substack{k \rightarrow \infty}} \biggm( \varrho(\mathcal{K}^{\ast}, \mathcal{K}^{(\ell_k)}) + \varepsilon_k \bigl\Vert \mathcal{K}^{\ast}- \mathcal{K}^{(\ell_k)} \bigr\Vert_{\prod\nolimits_{i=1}^n \mathscr{K}_i[\mathcal{X}, \mathcal{U}_i]}^2 \biggm) = 0, \label{Eq63}
\end{align}
and thereby provides $\mathcal{K}^{\ast} \in \prod\nolimits_{i=1}^n \mathscr{K}_i[\mathcal{X}, \mathcal{U}_i]$ is a fixed-point for the mapping $\varrho(\mathcal{K}^{\ast}, \mathcal{K})$, i.e., $\inf_{\substack {\mathcal{K}} \in \prod\nolimits_{i=1}^n \mathscr{K}_i[\mathcal{X}, \mathcal{U}_i]} \varrho(\mathcal{K}^{\ast}, \mathcal{K})=\varrho(\mathcal{K}^{\ast}, \mathcal{K}^{\ast})$, such that
\begin{align}
r_{i}(\mathcal{K}_i^{\ast}, \mathcal{K}_{\neg i}^{\ast}) \le r_{i}(\mathcal{K}_i, \mathcal{K}_{\neg i}^{\ast}), \,\, \forall \mathcal{K}_i \in \mathscr{K}_{i}[\mathcal{X}, \mathcal{U}_i],\,\, \forall i \in \{1, 2, \dots, n\}, \label{Eq64}
\end{align}
which shows that $\mathcal{K}^{\ast}=(\mathcal{K}_1^{\ast}, \mathcal{K}_2^{\ast}, \ldots, \mathcal{K}_n^{\ast}) \in \mathscr{K}$ is indeed a Nash equilibrium point for the noncooperative $n$-player principal eigenvalue game.\footnote{Such a map, whose fixed-point is an equilibrium, is called a Nash map for the game (see \cite{Glic52}).} With the admissible control strategies $u_i^{\ast}(\cdot) = \bigl(\mathcal{K}_i^{\ast} x^{\epsilon}\bigr)(\cdot)$, $\forall t \ge 0$, for each $i =1, 2, \ldots, n$, the maximum closed invariant set $\Lambda_D^{\mathcal{K}^{\ast}} \subset D \cup \partial D$ contains all closed invariant sets from the set $D \cup \partial D$ (under the action of the $n$-tuple of equilibrium linear feedback operators with respect to the unperturbed multi-channel dynamical system) i.e., $\Lambda_D^{(\mathcal{K}_i, \mathcal{K}_{\neg i}^{\ast})} \subset \Lambda_D^{\mathcal{K}^{\ast}}$, $\forall \mathcal{K}_i \in \mathscr{K}_{i}[\mathcal{X}, \mathcal{U}_i]$, $\forall i \in \{1,2, \ldots, n\}$ (see also Proposition~\ref{P1}). This completes the proof of Proposition~\ref{P3}. \qed

\end{document}